\documentclass[12pt,twoside]{amsart}
\usepackage{amssymb}
\usepackage{amscd}
\usepackage{hyperref}
\usepackage{color}

\usepackage[normalem]{ulem} 
\newcommand\redout{\bgroup\markoverwith
{\textcolor{red}{\rule[.5ex]{2pt}{0.4pt}}}\ULon}

\title[Mori fibre spaces in positive characteristic]
{Pathologies on Mori fibre spaces in positive characteristic} 
\author{Hiromu Tanaka} 
\subjclass[2010]{14E30.}
\keywords{minimal model program, Mori fibre space, positive characteristic}
\address{Graduate School of Mathematical Sciences, 
The University of Tokyo, 
3-8-1 Komaba, Meguro-ku, Tokyo 153-8914, JAPAN} 
\email{tanaka@ms.u-tokyo.ac.jp}

\newcommand{\red}[0]{{\operatorname{red}}}

\newcommand{\Proj}[0]{{\operatorname{Proj}}}
\newcommand{\Spec}[0]{{\operatorname{Spec}}}

\newcommand{\Supp}[0]{{\operatorname{Supp}}}

\newcommand{\Ex}[0]{{\operatorname{Ex}}}
\newtheorem{thm}{Theorem}[section]
\newtheorem{lem}[thm]{Lemma}

\newtheorem{prop}[thm]{Proposition}

\newtheorem*{claim}{Claim}  
\newtheorem{step}{Step}

\theoremstyle{definition}

\newtheorem{dfn}[thm]{Definition}

\newtheorem{rem}[thm]{Remark}
      
\newtheorem{nota}[thm]{Notation}

\makeatletter
  
  \@addtoreset{equation}{thm}
  \makeatother

\newcommand{\MO}{\mathcal{O}}
\newcommand{\R}{\mathbb{R}}
\newcommand{\Q}{\mathbb{Q}}
\newcommand{\Z}{\mathbb{Z}}
\newcommand{\F}{\mathbb{F}}

\newcommand{\q}{\mathfrak{q}}
\begin{document}

\maketitle

\begin{abstract}
We show that there exist Mori fibre spaces 
whose total spaces are klt but bases are not. 
We also construct Mori fibre spaces which have relatively non-trivial torsion line bundles. 
\end{abstract}

\tableofcontents

\section{Introduction}

Given an algebraic variety $X$, 
the minimal model conjecture implies that 
$X$ is birational to either a minimal model or a Mori fibre space. 
The purpose of this paper is 
to find some phenomena on Mori fibre spaces 
that occur only in positive characteristic. 
Originally, the advantage of Mori fibre spaces is their simple structure, 
which allows us to reduce some problems to the study of their fibres and bases. 
For instance, given a Mori fibre space $f:X \to S$ from a klt variety $X$ 
in characteristic zero, 
it is known that its base space $S$ is also klt 
(cf.  \cite[Theorem 0.2]{Amb05}, \cite[Corollary 4.6]{Fuj99}). 
Unfortunately, the same statement is no longer true in positive characteristic. 

\begin{thm}\label{intro-patho1}
Let $k$ be an algebraically closed field 
whose characteristic is two or three. 
Then there exists a projective $k$-morphism $f:V \to W$ of normal $k$-varieties that satisfies the following properties:
\begin{enumerate}
\item $V$ is a $4$-dimensional $\Q$-factorial klt variety over $k$, 
\item $W$ is a $3$-dimensional $\Q$-factorial variety over $k$ 
that is not klt, 
\item $f_*\MO_V=\MO_W$, $\rho(V/W)=1$, $-K_V$ is $f$-ample,
\item any fibre of $f$ is an irreducible scheme of dimension one, and 
there is a non-empty open subset $W^0$ of $W$ such that 
the fibre $V \times_W \Spec\,k(w)$ is isomorphic to $\mathbb P^1_{k(w)}$ 
for any point $w \in W^0$, where $k(w)$ denotes the residue field at $w$. 
\end{enumerate}
\end{thm}

A prominent property of Mori fibre spaces in characteristic zero 
is that any relatively numerically trivial Cartier divisor is trivial 
(cf. \cite[Lemma 3-2-5(2)]{KMM87}). 
We construct an example in positive characteristic that violates this property.

\begin{thm}\label{intro-patho2}
Let $k$ be an algebraically closed field 
whose characteristic $p$ is two or three. 
Then there exists a projective $k$-morphism $f:V \to W$ of normal $k$-varieties that satisfies the following properties:
\begin{enumerate}
\item $V$ is a $3$-dimensional $\Q$-factorial klt variety over $k$, 
\item $W$ is a smooth curve over $k$, 
\item $f_*\MO_V=\MO_W$, $\rho(V/W)=1$, $-K_V$ is $f$-ample, and 
\item there is a Cartier divisor $D$ on $V$ such that $D \not\sim_f 0$ and $pD \sim_f 0$. 
\end{enumerate}
\end{thm}

\begin{rem}
Since \cite[Lemma 3-2-5(2)]{KMM87} is 
a formal consequence of the relative Kawamata--Shokurov base point free theorem 
\cite[Theorem 3-1-1]{KMM87}, 
the same statement as in \cite[Theorem 3-1-1]{KMM87} 
does not hold in positive characteristic. 
\end{rem}

\subsection{Construction of examples}

Let us overview how to construct the examples appearing in 
Theorem~\ref{intro-patho1} and Theorem~\ref{intro-patho2}.

\subsubsection{Pathological surfaces over imperfect fields}

To find examples appearing in 
Theorem~\ref{intro-patho1} and Theorem~\ref{intro-patho2}, 
we first construct log del Pezzo surfaces over imperfect fields 
satisfying pathological properties as follows. 

\begin{thm}\label{intro-Fano-CY}
Let $k$ be an imperfect field whose characteristic $p$ is two or three. 
Then there exists a $k$-morphism $\rho:S \to C$ 
that satisfies the following properties:  
\begin{enumerate}
\item $S$ is a projective regular surface over $k$ and 
there is an effective $\Q$-divisor $\Delta_S$ 
such that $(S, \Delta_S)$ is klt and $-(K_S+\Delta_S)$ is ample, 
\item $C$ is a projective regular curve over $k$ with $K_C \sim 0$, 
\item $\rho$ is a $\mathbb P^1$-bundle, and  
\item there is a Cartier divisor $L$ on $C$ such that $L \not\sim 0$ and $pL \sim 0$. 
\end{enumerate}
\end{thm}

The surface $S$ in Theorem~\ref{intro-Fano-CY} 
is a log Fano variety dominating a Calabi--Yau variety. 
Such an example does not exist in characteristic zero 
(cf. \cite[Lemma 2.8]{PS09}, 
\cite[Theorem 5.1]{FG12}). 
For some related results in positive characteristic, we refer to \cite{Eji}. 

Let us overview the construction of 
$\rho:S \to C$ appearing in Theorem~\ref{intro-Fano-CY}. 
We take a regular cubic curve $C$  
that is not smooth and has a $k$-rational point $P$ around which 
$C$ is smooth over $k$. 
For example, if $k$ is the function field of a curve over an algebraically closed field, 
then $C$ is nothing but the generic fibre of a quasi-elliptic fibration 
equipped with a section. 
Since we have that $H^1(C, \MO_C(-P)) \neq 0$ by Serre duality, 
a nonzero element $\xi$ of $H^1(C, \MO_C(-P))$ 
induces a locally free sheaf $E$ of rank two. 
Then $S$ is the $\mathbb P^1$-bundle defined as $\mathbb P(E)$. 
In order to show that $S$ is log del Pezzo, 
one of the essential facts is 
that we can find a purely inseparable field extension 
$k \subset k'$ of degree $p$ such that $C \times_k k'$ is an integral but non-normal scheme and 
that its normalisation is isomorphic to $\mathbb P^1_{k'}$. 
Since the scheme-theoretic inverse image $\varphi^{-1}(P)$ of 
the $k$-rational point $P$ 
is a $k'$-rational point, we have that 
$$H^1(\mathbb P^1_{k'}, \MO_{\mathbb P^1_{k'}}(-\varphi^*P))
=H^1(\mathbb P^1_{k'}, \MO_{\mathbb P^1_{k'}}(-1))=0,$$
where $\varphi^*P$ denotes the pull-back of the Cartier divisor $P$. 
This implies that the pull-back $\varphi^*\xi$ is zero. 
This property plays a crucial role in our construction. 
For more details, see Section~\ref{s-dP}. 

\subsubsection{Proofs of the theorems}

Let us overview some of the ideas of the proofs of 
Theorem~\ref{intro-patho1} and Theorem~\ref{intro-patho2}. 

First let us treat the latter one: Theorem~\ref{intro-patho2}. 
This is a consequence of Theorem~\ref{intro-Fano-CY}. 
Indeed, for an algebraically closed field $k$ 
whose characteristic is two or three, 
it follows from Theorem~\ref{intro-Fano-CY} 
that we get a log del Pezzo surface $(S, \Delta_S)$ over $k(t)$ 
which has a non-trivial $p$-torsion Cartier divisor. 
Then we can spread it out 
over some non-empty open subset $W$ of $\mathbb A^1_k$, 
i.e. there is a morphism $V \to W$ with $V \times_W \Spec\,k(t)=S$. 
Although this example does not satisfy the property $\rho(V/W)=1$, 
we may assume this condition by contracting an appropriate curve on $S$ in advance. 
For more details, see Subsection~\ref{ss-MFS1}. 

Second, let us overview the proof of Theorem~\ref{intro-patho1}. 
To this end, we first find a similar example over imperfect fields 
(cf. Theorem~\ref{t-imperfect-nonklt}). 
The basic idea is to take cones over $\rho:S \to C$. 
However, there is no morphism between cones. 
What we will actually do is to take $\mathbb P^1$-bundles functorially 
for an ample divisor $M_C$ on $C$: 
$$X:=\mathbb P_S(\MO \oplus \MO(\rho^*M_C)) \to 
\mathbb P_C(\MO \oplus \MO(M_C))=:W_0.$$
Let $W_0 \to W$ 
be the birational contraction of the section $C^-$ of $W_0 \to C$ 
with negative self-intersection number. 
Since $K_C \sim 0$, $W$ is not klt. 

If there was a divisorial contraction whose exceptional locus is 
the pull-back of $C^-$, 
then the resulting variety would be what we are looking for.  
Although we can not hope this, 
we will get close to this situation by running a suitable minimal model program. 
To this end, 
we first construct a minimal model program after taking a purely inseparable cover of $X$, 
and descend it to $X$ after that. 
For more details, see Subsection~\ref{ss-MFS2}.

\medskip

\textbf{Acknowledgement:} 
The author would like to thank P. Cascini, S. Ejiri, Y. Gongyo and J. Witaszek 
for useful comments and answering questions. 
He also thanks the referee for many constructive suggestions and 
reading the manuscript carefully. 
The author was funded by EPSRC and 
the Grant-in-Aid for Scientific Research (KAKENHI No. 18K13386).

\section{Preliminaries}

\subsection{Notation}\label{ss-notation}

In this subsection, we summarise the notation used in this paper. 

\begin{itemize}
\item We will freely use the notation and terminology in \cite{Har77} 
and \cite{Kol13}. 
\item For a scheme $X$, its {\em reduced structure} $X_{\red}$ 
is the reduced closed subscheme of $X$ such that the induced morphism 
$X_{\red} \to X$ is surjective. 
\item For an integral scheme $X$, 
we define the {\em function field} $K(X)$ of $X$ 
as $\MO_{X, \xi}$ for the generic point $\xi$ of $X$. 
\item For a field $k$, 
$X$ is a {\em variety over} $k$ or a $k$-{\em variety} if 
$X$ is an integral scheme that is separated and of finite type over $k$. 
We say that $X$ is a {\em curve} over $k$ or a $k$-{\em curve} 
(resp. a {\em surface} over $k$ or a $k$-{\em surface}, 
resp. a {\em threefold} over $k$) 
if $X$ is a $k$-variety of dimension one (resp. two, resp. three). 
\item We say that two schemes $X$ and $Y$ over a field $k$ are $k$-{\em isomorphic} 
if there exists an isomorphism $\theta\colon X \to Y$ of schemes 
such that both $\theta$ and $\theta^{-1}$ commute with the structure morphisms: 
$X \to \Spec\,k$ and $Y \to \Spec\,k$. 
\end{itemize}

\begin{dfn}\label{def-intersection}
Let $k$ be a field. 
\begin{enumerate}
\item{Let $C$ be a proper curve over $k$. 
Let $M$ be an invertible sheaf on $C$. 
It is well-known that 
$$\chi(C, mM)=\dim_k(H^0(C, mM))-\dim_k(H^1(C, mM)) \in\Z[m]$$ 
and that the degree of this polynomial is at most one 
(cf. \cite[Ch I, Section 1, Theorem on page 295]{Kle66}). 
We define the {\em degree} of $M$ over $k$, 
denoted by $\deg_k M$ or $\deg M$, as the coefficient of $m$. }
\item{Let $X$ be a separated scheme of finite type over $k$, 
let $L$ be an invertible sheaf on $X$, and 
let $C \hookrightarrow X$ be a closed immersion over $k$ 
from a proper $k$-curve $C$. 
We define the {\em intersection number} over $k$, 
denoted by $L\cdot_k C$ or $L \cdot C$, as $\deg_{k}(L|_C).$}
\end{enumerate}
\end{dfn}

\subsection{Properties spreading out from the generic fibre}

In this subsection, we summarise 
some properties extendable from the generic fibre: 
Lemma~\ref{l-rel-properties}. 
Also, we give a criterion 
for projective morphisms of generically relative Picard number one (Lemma \ref{l-picard-one}). 
To this end, we establish two auxiliary lemmas: 
Lemma~\ref{l-Pn-bundle} and Lemma~\ref{l-rel-inv}.

\begin{lem}\label{l-rel-properties}
Let $k$ be a field. 
Let $f:X \to Y$ be a projective $k$-morphism of normal $k$-varieties 
with $f_*\MO_X=\MO_Y$. 
\begin{enumerate}
\item 
Assume that $k$ is algebraically closed. 
Then the generic fibre $X_{K(Y)}$ is $\Q$-factorial if and only if 
there is a non-empty open subset $Y'$ of $Y$ such that $X \times_Y Y'$ is $\Q$-factorial. 
\item 
Let $\Delta$ be an effective $\Q$-divisor 
such that $K_X+\Delta$ is $\Q$-Cartier. 
Assume that there is a log resolution of $(X_{K(Y)}, \Delta|_{X_{K(Y)}})$. 
Then $(X_{K(Y)}, \Delta|_{X_{K(Y)}})$ is klt (resp. log canonical) 
if and only if 
there is a non-empty open subset $Y'$ of $Y$ such that $(X \times_Y Y', \Delta|_{X \times_Y Y'})$ is klt (resp. log canonical). 
\end{enumerate}
\end{lem}

\begin{proof}
The assertion (1) holds by \cite[the third Theorem in Introduction]{BGS}. 
We now show (2). 
After shrinking $Y$, we may assume that there is a log resolution $g:Z \to X$ of $(X, \Delta)$. 
We define a $\Q$-divisor $\Delta_Z$ on $Z$ by $K_Z+\Delta_Z=g^*(K_X+\Delta)$. 
Then \cite[Corollary 2.13(2)]{Kol13} implies the following:  
\begin{enumerate}
\item[(a)] $(X, \Delta)$ is klt if and only if all the coefficients of $\Delta_Z$ are less than $1$. 
\item[(b)] $(X_{K(Y)}, \Delta|_{X_{K(Y)}})$ is klt if and only if 
all the coefficients of $\Delta_Z|_{Z_{K(Y)}}$ are less than $1$. 
\end{enumerate}
Shrinking $Y$ again, we may assume that any irreducible component of $\Delta_Z$ dominates $Y$. 
Then all the coefficients of $\Delta_Z$ are less than $1$ if and only if 
all the coefficients of $\Delta_Z|_{Z_{K(Y)}}$ are less than $1$. 
Thanks to (a) and (b), $(X, \Delta)$ is klt if and only if $(X_{K(Y)}, \Delta|_{X_{K(Y)}})$ is klt. 
\end{proof}

\begin{lem}\label{l-Pn-bundle}
Let $k$ be a field. 
Let $f:X \to Y$ be a projective $k$-morphism of normal $k$-varieties. 
Assume that the generic fibre $X_{K(Y)}$ is $K(Y)$-isomorphic to $\mathbb P^n_{K(Y)}$ for some non-negative integer $n$. 
Then there exists a non-empty open subset $Y'$ of $Y$ such that 
the fibre $X_y$ is $k(y)$-isomorphic to $\mathbb P^n_{k(y)}$ for any point $y \in Y'$. 
\end{lem}

\begin{proof} 
Replacing $Y$ by a non-empty open subset, we may assume that the following properties hold: 
\begin{enumerate}
\item $f$ is a smooth morphism and $f_*\MO_X=\MO_Y$, 
\item $-K_X$ is $f$-ample,  
\item the tangent bundle $T_{X_y}$ is ample for any $y \in Y$ 
(cf. \cite[Proposition 4.4]{Har66}, \cite[Proposition 6.1.9]{Laz04}), and 
\item there is a section of $f$, i.e. there exists a closed immersion $j:Y_1 \to X$ such that 
the composite morphism $Y_1 \to X \to Y$ is an isomorphism. 
\end{enumerate}
Fix $y \in Y$ and let $X_{\overline{k(y)}}$ be the base change of the fibre $X_y$ to its algebraic closure $\overline{k(y)}$. 
Since $X_{\overline{k(y)}}$ is a smooth projective variety whose tangent bundle 
$T_{X_{\overline{k(y)}}}$ is ample, 
we have that $X_{\overline{k(y)}}$ is $\overline{k(y)}$-isomorphic to 
$\mathbb P^n_{\overline{k(y)}}$ by \cite[Theorem 8]{Mor79}. 
This implies that $X_y$ is a Severi--Brauer variety. 
Since $X_y$ has a $k(y)$-rational point by (4), 
we have that $X_y$ is $k(y)$-isomorphic to $\mathbb P^n_{k(y)}$ by \cite[Theorem 5.1.3]{GS06}. 
\end{proof}

\begin{lem}\label{l-rel-inv}
Let $k$ be a field of characteristic $p>0$. 
Consider a commutative diagram of projective $k$-morphisms of normal 
$k$-varieties 
$$\begin{CD}
X' @>\alpha >> X\\
@VVf'V @VVfV\\
Y' @>\beta >> Y,
\end{CD}$$
where $\alpha$ and $\beta$ are finite universal homeomorphisms. 
Then $\rho(X/Y)=\rho(X'/Y')$. 
\end{lem}

\begin{proof}
As $X$ and $Y$ are normal, 
the $e$-th iterated absolute Frobenius morphisms $F^e_X:X \to X$ and $F^e_Y:Y \to Y$ 
can be considered as the normalisation of $X$ in $K(X)^{1/p^e}$ and $Y$ in $K(Y)^{1/p^e}$, respectively. 
Since there is a positive integer $e$ such that $K(X)^{1/p^e} \supset K(X') \supset K(X)$ 
and $K(Y^{1/p^e}) \supset K(Y') \supset K(Y)$,  
the $e$-th iterated absolute Frobenius morphisms 
$F^e_X:X \to X$ and $F^e_Y:Y \to Y$ 
factor through $\alpha$ and $\beta$ respectively: 
$$F^e_X:X \xrightarrow{\widetilde \alpha} X' \xrightarrow{\alpha} X, \quad F^e_Y:Y \xrightarrow{\widetilde \beta} Y' \xrightarrow{\beta} Y.$$
Then we have that $\rho(X/Y) \leq \rho(X'/Y')$. 
The opposite inequality follows from the fact that 
the $e$-th iterated absolute Frobenius morphisms 
$F^{e}_{X'}:X' \to X'$ and $F^e_{Y'}:Y' \to Y'$ factor through $\widetilde \alpha$ and 
$\widetilde \beta$, respectively. 
\end{proof}

\begin{lem}\label{l-picard-one}
Let $k$ be a field of characteristic $p>0$. 
Let $f:X \to Y$ be a projective $k$-morphism of normal $k$-varieties 
with $f_*\MO_X=\MO_Y$. 
Assume that there exists a finite universal homeomorphism $\varphi:\mathbb P_L^n \to X_{K(Y)}$ 
over $K(Y)$ for some finite purely inseparable extension $K(Y) \subset L$. 
Then there exists a non-empty open subset $Y'$ of $Y$ such that $\rho(X'/Y')=1$ for $X':=X \times_Y Y'$. 
\end{lem}

\begin{proof}
We have a commutative diagram 
$$\begin{CD}
X_1:=\mathbb P^n_L @>\varphi>> X_{K(Y)} @>>> X\\
@VVf_1 V @VVf_{K(Y)}V @VVfV\\
Y_1:=\Spec\,L @>>> \Spec\,K(Y) @>>> Y,
\end{CD}$$
where the right square is cartesian. 
All the schemes in the left square are projective over $K(Y)$, hence 
after shrinking $Y$, we can find a commutative diagram of 
projective $k$-morphisms of normal 
$k$-varieties 
$$\begin{CD}
X_2 @>\alpha >> X\\
@VVf_2 V @VVfV\\
Y_2 @>\beta >> Y,
\end{CD}$$
where $\alpha$ and $\beta$ are finite universal homeomorphisms, 
$K(Y_2) =L$ and the generic fibre of $f_2$ is $K(Y_2)$-isomorphic to 
$\mathbb P_{K(Y_2)}^n$. 
Then the assertion follows 
from Lemma \ref{l-Pn-bundle} and Lemma \ref{l-rel-inv}. 
\end{proof}

\subsection{Varieties of Fano type}

In this subsection, we recall the definition of varieties of Fano type and 
one of basic properties (Lemma~\ref{l-fano-birat}). 

\begin{dfn}
Let $k$ be a field. 
A projective normal $k$-variety $X$ is {\em of Fano type} 
if there is an effective $\Q$-divisor $\Delta$ such that 
$(X, \Delta)$ is klt and $-(K_X+\Delta)$ is ample. 
In this case, $(X, \Delta)$ is called {\em log Fano}. 
We say that $(X, \Delta)$ is {\em log del Pezzo} if $X$ is log Fano and $\dim X=2$. 
\end{dfn}

\begin{lem}\label{l-klt-perturb}
Let $k$ be a field of characteristic $p>0$ such that $[k:k^p]<\infty$. 
Let $(X, \Delta)$ be a projective klt pair over $k$ such that $\dim X\leq 3$. 
If $A$ is a nef and big $\Q$-Cartier $\Q$-divisor on $X$, 
then there exists an effective $\Q$-Cartier $\Q$-divisor $A'$ such that 
$A \sim_{\Q} A'$ and $(X, \Delta+A)$ is klt. 
\end{lem}

\begin{proof}
Thanks to the assumptions $[k:k^p]<\infty$ and $\dim X \leq 3$, 
we may freely use log resolutions by \cite{CP08, CP09}. 
Then we can apply the same argument as in \cite[Lemma 2.8]{GNT}. 
\end{proof}

\begin{lem}\label{l-fano-birat}
Let $k$ be a field of characteristic $p>0$. 
Let $X$ and $Y$ be projective normal varieties over $k$. 
Assume that a rational map $f:X \dashrightarrow Y$ over $k$ 
satisfies one of the following properties. 
\begin{enumerate}
\item $f$ is a birational morphism. 
\item $f$ is a birational map which is an isomorphism in codimension one. 
\end{enumerate}
If $X$ is of Fano type, $[k:k^p]<\infty$ and $\dim X \leq 3$, 
then $Y$ is of Fano type. 
\end{lem}

\begin{proof}
For both the cases, we can apply the same argument as in 
\cite[Lemma 2.4]{Bir16}. 
However, we give a proof only for the case (1) since our setting 
differs from \cite[Lemma 2.4]{Bir16}. 
Thanks to the assumptions $[k:k^p]<\infty$ and $\dim X \leq 3$, 
we may freely use log resolutions by \cite{CP08, CP09}. 

Since $X$ is of Fano type, 
we can find an effective $\Q$-divisor $\Delta$ on $X$ such that 
$(X, \Delta)$ is klt and $-(K_X+\Delta)$ is ample. 
By Lemma \ref{l-klt-perturb}, 
we can find an effective $\Q$-divisor $A_X$ on $X$ 
such that $-(K_X+\Delta) \sim_{\Q} A_X$ and $(X, \Delta+A_X)$ is klt. 
Taking the push-forward by $f$, we have that 
$$-(K_Y+f_*\Delta+f_*A_X) \sim_{\Q} 0.$$
Then it holds that 
\[
K_X+\Delta+A_X=f^*(K_Y+f_*\Delta+f_*A_X),
\]
hence $(Y, f_*\Delta+f_*A_X)$ is klt. 
Since $f_*A_X$ is big, we can write 
$$f_*A_X=A_Y+E$$
for some ample $\Q$-Cartier $\Q$-divisor $A_Y$ and an effective $\Q$-divisor $E$. 
Therefore, it holds that 
$$-(K_Y+f_*\Delta+(1-\epsilon)f_*A_X+\epsilon E) \sim_{\Q} \epsilon A_Y$$
for any rational number $\epsilon$. 
Since log resolutions exist, the pair 
$$(Y, f_*\Delta+(1-\epsilon)f_*A_X+\epsilon E)$$
is klt if $\epsilon$ is a sufficiently small rational number. 
Hence, $Y$ is of Fano type. 
\end{proof}

\begin{rem}
The assumptions 
$[k:k^p]<\infty$ and $\dim X \leq 3$ in Lemma~\ref{l-fano-birat} is used 
only to assure the existence of log resolutions \cite{CP08, CP09}. 
\end{rem}

\subsection{Jacobian criterion for regularity}

For later use, 
we summarise results for regularity of some explicit varieties 
that follow from the Jacobian criterion. 

\begin{lem}\label{l-Jacobian}
Let $k$ be a field of characteristic $p>0$. 
Take elements $s, t \in k \setminus k^p$.  
Then the following hold. 
\begin{enumerate}
\item If $p=2$, 
then $\Spec\,k[x, y]/(x^2+ty^2)$ is regular outside the origin $\{(0, 0)\}$. 
\item If $p=2$, then $k[x, y]/(tx^2+1)$ is regular. 
\item If $p=2$ and $[k(s^{1/2}, t^{1/2}):k]=4$, then $k[x, y]/(sx^2+ty^2+1)$ is regular. 
\item If $p=2$, then $\Proj\,k[x, y, z]/(y^2z+x^3+sxz^2)$ is regular. 
\item If $p=3$, then $\Proj\,k[x, y, z]/(-y^2z+x^3+sz^3)$ is regular.
\end{enumerate}
\end{lem}

\begin{proof}
Since all the proofs are quite similar, we only prove (1). 
We consider the following open subset of $\Spec\,k[x, y]/(x^2+ty^2)$: 
$$D(y)=\Spec\,k[x, y, z]/(x^2+ty^2, zy+1) \simeq 
\Spec\,k[x, y, y^{-1}]/(x^2y^{-2}+t).$$
We can find an $\F_2$-derivation $D_1$ of $k[x, y, z]$ with $D_1(t)=1$ 
by $t \not\in k^2$ and \cite[Theorem 26.5]{Mat89}. 
We have that the ring $k[x, y, z]/(x^2+ty^2, zy+1)$ is regular 
by applying the Jacobian criterion \cite[Proposition 22.6.7(iii)]{EGAIV} 
for $k_0:=\F_2$, $B=k[x, y, z]$, 
$\q$ is a prime ideal of $B$ containing $(x^2+ty^2, zy+1)$, 
$f_1:=x^2+ty^2$, and the $\F_2$-derivation $D_1$ defined above. 
It follows from the same argument that 
also the open subset $D(x)$ of $\Spec\,k[x, y]/(x^2+ty^2)$ is regular. 
To summarise, $\Spec\,k[x, y]/(x^2+ty^2)$ is regular outside the origin $(0, 0)$. 
\end{proof}

\subsection{Slc-ness of conics}

The purpose of this subsection is to show that 
any plane conic curve is semi log canonical (Lemma~\ref{l-conic}). 
We start with a typical case in characteristic two. 

\begin{lem}\label{e-conic/p=2}
Let $k$ be a field of characteristic two 
with an element $t \in k$ such that $t \not\in k^2$. 
Let 
$$Z:=\Spec\,k[x, y]/(x^2+ty^2).$$
Then $Z$ is a semi log canonical curve. 
\end{lem}

\begin{proof}
By Lemma \ref{l-Jacobian}(1), $Z$ is regular outside the origin. 
By using the assumption: $t \not\in k^2$, 
we can directly check that $Z$ has a node at the origin (cf. \cite[1.41]{Kol13}). 
It follows from \cite[Proposition 3.6]{Tan16} that $Z$ is semi log canonical. 
\end{proof}

\begin{lem}\label{l-conic}
Let $k$ be a field and let $Z=\Spec\,k[x, y]/(f)$, 
where $Z$ is reduced and $f$ is a polynomial of degree two. 
Then $Z$ is semi log canonical. 
\end{lem}

\begin{proof}
We only treat the case where the characteristic of $k$ is two, 
since otherwise the problem is easier. 
We may assume that $k$ is separably closed. 
We can write 
$$f=a_{20}x^2+a_{11}xy+a_{02}y^2+a_{10}x+a_{01}y+a_{00}$$ 
where $a_{ij} \in k$. 

\setcounter{step}{0}

\begin{step}\label{s1-conic}
If $a_{11} \neq 0$, then $Z$ is semi log canonical. 
\end{step}

\begin{proof}(of Step \ref{s1-conic}) 
Assume that $a_{11} \neq 0$. 
Since $k$ is separably closed, we can write 
$a_{20}x^2+a_{11}xy+a_{02}y^2=l_1l_2$ for 
some homogeneous polynomials $l_1$ and $l_2$ of degree one 
with $(l_1, l_2)=(x, y)$. 
In particular, we may assume that $a_{20}=a_{02}=0$. 
After applying a suitable linear transform, we may assume that $f=xy+b$ for some $b \in k$. 
Then $Z$ is semi log canonical. 
This completes the proof of Step \ref{s1-conic}. 
\end{proof}

\begin{step}\label{s2-conic}
If $a_{11} = 0$ and $a_{10} \neq 0$, then $Z$ is smooth over $k$. 
\end{step}

\begin{proof}(of Step \ref{s2-conic}) 
We first prove that we may assume that $a_{00}=0$. 
If $a_{20}=0$, then we are reduced to the case where $a_{00}=0$, 
after applying the linear transform: $a_{10}x+a_{00} \mapsto x$, $y \mapsto y$. 
Assume that $a_{20} \neq 0$. Since $k$ is separably closed, 
$a_{20}x^2+a_{10}x+a_{00}=0$ has two distinct solutions if $a_{20} \neq 0$. 
In particular, applying the linear transform $x+\alpha \mapsto x$, $y \mapsto y$ 
for a solution $\alpha \in k$, 
we may assume that $a_{00}=0$. 

\medskip

We assume that $a_{00}=0$. Thus, we can write 
$$f=a_{20}x^2+a_{02}y^2+a_{10}x+a_{01}y.$$ 
Since $a_{10} \neq 0$, we have that $Z \simeq k[x, y]/(x+bx^2+cy^2)$ 
for some $b, c \in k$. 
It follows from the Jacobian criterion for smoothness 
that $Z$ is smooth.  
This completes the proof of Step \ref{s2-conic}. 
\end{proof}

\begin{step}\label{s3-conic}
If $a_{11} = a_{10}=a_{01}=0$, then $Z$ is semi log canonical. 
\end{step}

\begin{proof}(of Step \ref{s3-conic}) 
If $a_{00}=0$, then we may assume that $a_{20}=1$ by symmetry, i.e. 
$f=x^2+a_{02}y^2$. 
Since $Z$ is reduced, we have that $a_{02} \not\in k^2$. 
Thus $Z$ is semi log canonical by Lemma~\ref{e-conic/p=2}. 

Thus we may assume that $a_{00} \neq 0$. 
Replacing $f$ by $f/{a_{00}}$ we get 
$$f=a_{20}x^2+a_{02}y^2+1.$$
Since $f$ is reduced, we may assume that $a_{20} \not\in k^2$, 
hence $[k(a_{20}^{1/2}):k]=2$. 
This implies that 
$[k(a_{20}^{1/2}, a_{02}^{1/2}):k]$ is equal to either $4$ or $2$. 
If $[k(a_{20}^{1/2}, a_{02}^{1/2}):k]=4$, 
then $Z$ is regular by Lemma~\ref{l-Jacobian}(3). 

Thus we may assume that $[k(a_{20}^{1/2}, a_{02}^{1/2}):k]=2$. 
We have that $a_{02} \in k^2(a_{20})=k^2 \oplus k^2a_{20}$, 
hence $a_{02}=b^2+c^2a_{20}$ for some $b, c \in k$. 
We get  
$$f=a_{20}x^2+a_{02}y^2+1=a_{20}x^2+(b^2+c^2a_{20})y^2+1=a_{20}(x+cy)^2+(by+1)^2.$$
After replacing $x+cy$ by $x$, we can write 
$$f=a_{20}x^2+(by+1)^2.$$
If $b=0$, then $Z$ is regular by Lemma~\ref{l-Jacobian}(2). 
If $b \neq 0$, then we get 
$k[x, y]/(f) \simeq k[x, y]/(a_{20}x^2+y^2)$. 
It follows from Lemma~\ref{e-conic/p=2} that $Z$ is semi log canonical. 
This completes the proof of Step \ref{s3-conic}. 
\end{proof}

Step \ref{s1-conic}, Step \ref{s2-conic} and Step \ref{s3-conic} complete 
the proof of Lemma \ref{l-conic}. 
\end{proof}





\section{Pathological surfaces over imperfect fields}\label{s-dP}

\subsection{Construction in a general setting}\label{ss-construction}

In this subsection, we give a criterion to find a log Fano variety 
dominating a Calabi--Yau variety (Proposition~\ref{p-criterion}). 
Although our construction is analogous to a standard one 
over an algebraically closed field (cf. \cite{Tan72}, \cite{Muk13}), 
we give details of them because our setting is more general and 
our base field is not necessarily algebraically closed.

\begin{nota}\label{n-general}
Let $k$ be a field of characteristic $p>0$. 
Assume that there exist a $k$-morphism $\varphi:C' \to C$ of regular projective $k$-varieties 
and a Cartier divisor $D$ on $C$ which satisfy the following properties. 
\begin{enumerate}
\item $\varphi$ is a finite universal homeomorphism of degree $p$. 
\item There is a nonzero element $\xi \in H^1(C, \MO_C(D))$ whose pull-back 
$\varphi^*(\xi) \in H^1(C', \MO_{C'}(\varphi^*D))$ is zero. 
\end{enumerate}
\end{nota}

The element $\xi$ induces a locally free sheaf $E$ of rank two on $C$ 
equipped with the following exact sequence 
that does not split:  
\begin{equation}\label{ex}
0 \to \MO_C(D) \xrightarrow{\alpha} E \xrightarrow{\beta} \MO_C \to 0.
\end{equation}
By our assumption (2), the pull-back of this sequence to $C'$ splits: 
$\varphi^*E \simeq \MO_{C'} \oplus \MO_{C'}(\varphi^*D)$. 
We set 
$$S:=\mathbb P_C(E),\quad S':=\mathbb P_{C'}(\varphi^*E)$$
and obtain a cartesian diagram: 
$$\begin{CD}
S' @>\psi >> S\\
@VV\rho' V @VV\rho V\\
C' @>\varphi>> C.
\end{CD}$$

The surjection $\beta:E \to \MO_C$ in (\ref{ex}) induces a section $C_1$ 
of $\rho$. 
We set $C'_1:=\psi^*C_1$ which is a section of $\rho'$. 
We have another section $C'_2$ of $\rho'$ 
corresponding to the surjection: 
$$\varphi^*E \simeq \MO_{C'} \oplus \MO_{C'}(\varphi^*D) \to \MO_{C'}(\varphi^*D),$$ 
where the latter homomorphism is the natural projection. 
We set $C_2$ to be the reduced closed subscheme of $S$ that is set-theoretically equal to $\psi(C'_2)$. 

\begin{lem}\label{l-p-cover}
We use Notation \ref{n-general}. 
Then $C_2$ is an integral scheme and the induced morphism $\rho_{C_2}:C_2 \to C$ is a finite universal homeomorphism of degree $p$. 
\end{lem}

\begin{proof}
Since $\psi$ is a universal homeomorphism, 
we have that $C_2$ is an integral scheme. 
Since the induced composite morphism 
$$C'_2 \xrightarrow{\psi} C_2 \xrightarrow{\rho_{C_2}} C$$ 
is a finite universal homeomorphism of degree $p$, 
the degree of the latter morphism $\rho_{C_2}:C_2 \to C$ is equal to either $1$ or $p$. 
It suffices to show that the latter case actually happens. 

Assuming that $\rho_{C_2}:C_2 \to C$ is of degree one i.e. birational, let us derive a contradiction. 
Since $\rho_{C_2}:C_2 \to C$ is a finite birational morphism and $C$ is normal, 
it follows that $\rho_{C_2}$ is an isomorphism. 
Thus $C_2$ is a section of $\rho$, hence 
it is corresponding to a surjective $\MO_C$-module homomorphism 
$$\gamma:E \to \MO_C(\widetilde D)$$
for some Cartier divisor $\widetilde D$ on $C$. 
Since $C'_2=C_2 \times_S S'$, we have that 
the pull-back $\varphi^*(\gamma \circ \alpha)$ is an isomorphism. 
It follows from the faithfully flatness of $\varphi$ that $\gamma \circ \alpha$ is an isomorphism. 
This implies that the sequence (\ref{ex}) splits, which is a contradiction. 
\end{proof}

\begin{lem}\label{l-pm-sections}
We use Notation \ref{n-general}. 
Then the following hold. 
\begin{enumerate}
\item $\MO_{S'}(C'_1)|_{C'_1} \simeq (\rho'_{C'_1})^*\MO_{C'}(-\varphi^*D)$, 
where $\rho'_{C'_1}:C'_1 \to C'$ is the induced morphism. 
\item $\MO_{S'}(C'_2)|_{C'_2} \simeq (\rho'_{C'_2})^*\MO_{C'}(\varphi^*D)$, 
where $\rho'_{C'_2}:C'_2 \to C'$ is the induced morphism. 
\end{enumerate}
\end{lem}

\begin{proof}
We can apply the same argument as in \cite[Ch. V, Proposition 2.6]{Har77}. 
\end{proof}

\begin{lem}\label{l-cano-div}
We use Notation \ref{n-general}. 
Then the following $\Q$-linear equivalence holds: 
$$-K_S \sim_{\Q} \frac{2}{p} C_2+\rho^*(-D-K_C).$$
\end{lem}

\begin{proof}
Since the induced morphism $\rho_{C_2}:C_2 \to C$ is of degree $p$ 
by Lemma~\ref{l-p-cover}, we have that 
\begin{equation}\label{e-cano-div1}
K_{S/C}+\frac{2}{p} C_2 \sim_{\Q} \rho^*(L)
\end{equation}
for some $\Q$-divisor $L$ on $C$. 
Taking the pull-back $\psi^*$ of (\ref{e-cano-div1}), we get 
\begin{equation}\label{e-cano-div2}
K_{S'/C'}+\frac{2}{p} \psi^*C_2 \sim_{\Q} \rho'^*\varphi^*(L).
\end{equation}
Since $K_{S'}+C'_1+C'_2+ \sim \rho'^*K_{C'}$, we have that 
\begin{equation}\label{e-cano-div3}
-C'_1-C'_2+\frac{2}{p} \psi^*C_2 \sim_{\Q} \rho'^*\varphi^*(L).
\end{equation}
It holds that 
$$(\rho'_{C'_1})^*\varphi^*D \sim_{\Q} 
-C'_1|_{C'_1} \sim_{\Q} (\rho'_{C'_1})^*\varphi^*(L),$$
where the first $\Q$-linear equivalence follows from Lemma \ref{l-pm-sections}(1) 
and the second one is obtained by restricting (\ref{e-cano-div3}) to $C'_1$. 
Since $\rho'_{C'_1}:C'_1 \to C'$ is an isomorphism 
and the absolute Frobenius morphism $F_C:C \to C$ factors through $\varphi:C' \to C$, 
we get the $\Q$-linear equivalence 
$pD=F_C^*(D) \sim_{\Q} F_C^*(L)=pL$, which implies  
\begin{equation}\label{e-cano-div4}
D \sim_{\Q} L.
\end{equation}
Substituting (\ref{e-cano-div4}) in (\ref{e-cano-div1}), 
we get  
$$-K_S\sim_{\Q} \frac{2}{p} C_2+\rho^*(-L-K_C) \sim_{\Q} \frac{2}{p} C_2+\rho^*(-D-K_C),$$
as desired. 
\end{proof}

\begin{prop}\label{p-criterion}
We use Notation \ref{n-general}. 
If $-D-K_C$ is ample and $(S, \frac{2}{p}C_2)$ is log canonical, then 
there exists an effective $\Q$-divisor $\Delta$ on $S$ such that $(S, \Delta)$ is klt and 
$-(K_S+\Delta)$ is ample. 
\end{prop}

\begin{proof}
Since $-D-K_C$ is ample, so is 
$$-\left(K_S+\middle(\frac{2}{p}-\epsilon \middle)C_2\right) \sim_{\Q} \epsilon C_2+\rho^*(-D-K_C)$$ 
for some rational number $\epsilon$ with $0 < \epsilon<\frac{2}{p}$, 
where the $\Q$-linear equivalence follows from Lemma~\ref{l-cano-div}. 
Set $\Delta:=(\frac{2}{p}-\epsilon)C_2$. 
Since $S$ is regular and  $(S, \frac{2}{p}C_2)$ is log canonical, 
we have that $(S, \Delta)$ is klt. 
\end{proof}

\subsection{Non-smooth K-trivial curves}

We summarise the properties of $K$-trivial curves that we will need later.

\begin{prop}\label{p-cubic}
Let $k$ be an imperfect field whose characteristic $p$ is two or three. 
Then there exists a projective regular curve $C$ over $k$ 
that satisfies the following properties:  
\begin{enumerate}
\item $K_C \sim 0$, 
\item the number of the $k$-rational points of $C$ is at least three, 
\item there is a purely inseparable field extension $k \subset k'$ of degree $p$ 
such that $C \times_k k'$ is an integral scheme which has 
a unique non-regular point $Q$, 
\item the normalisation $C'$ of $C \times_k k'$ 
is $k'$-isomorphic to $\mathbb P^1_{k'}$, and 
\item there is a Cartier divisor $L$ on $C$ such that $L \not\sim 0$ and $pL \sim 0$. 
\end{enumerate}
\end{prop}

\begin{proof}
Since $k$ is not perfect, we can find an element $t \in k$ with $t \not\in k^p$. 

First we treat the case where $p=2$. 
Consider the following equation, 
which is taken from \cite[Table 1 in page 243]{Ito94}: 
$$C:=\Proj\,k[x, y, z]/(y^2z+x^3+(t^3+t)xz^2).$$
We have that $C$ is regular by Lemma~\ref{l-Jacobian}(4). 
By the adjunction formula, (1) holds. 
The property (2) holds since all of $[0:1:0]$, $[0:0:1]$ and $[t+1:t^2+1:1]$ 
are $k$-rational points on $C$. 
Let $k':=k(\sqrt{t^3+t})$. 
The Jacobian criterion for smoothness implies that 
$C \times_k k'$ is smooth over $k'$ away from $Q:=[\sqrt{t^3+t}:0:1]$. 
Thus (3) holds. 
We can check that for the open set 
$$D_+(z)=\Spec\,k'[x, y]/(y^2+x^3+(t^3+t)x)$$
of $C\times_k k'$, 
its normalisation is isomorphic to $\mathbb A^1_{k'}$. 
Indeed, the integral closure of 
$k'[x, y]/(y^2+x^3+(t^3+t)x)=k'[\overline x, \overline y]$ is 
equal to $k'[\overline y/(\overline x+\sqrt{t^3+t})]$, where 
$\overline x$ and $\overline y$ are the images of $x$ and $y$, respectively. 
Thus (4) holds. 
The property (5) holds by setting $L:=P_1-P_2$, 
where $P_1$ and $P_2$ are $k$-rational points around which $C$ is smooth.

Second we assume that $p=3$. 
Let 
$$C:=\Proj\,k[x, y, z]/(-y^2z+x^3+t^2z^3).$$
All of $[0:1:0]$, $[0:t:1]$ and $[0:-t:1]$ 
are $k$-rational points on $C$. 
By Lemma~\ref{l-Jacobian}(5), $C$ is regular. 
We omit the remaining proof since it is similar to but easier than 
the one for the case where $p=2$. 
\end{proof}

\subsection{Log del Pezzo surfaces}

\begin{nota}\label{n-dP}
Let $k$ be an imperfect field whose characteristic $p$ is two or three. 
Let $C$ be a projective regular curve over $k$ as in Proposition~\ref{p-cubic}. 
By Proposition~\ref{p-cubic}(3)(4), 
there is a purely inseparable field extension 
$k \subset k'$ of degree $p$ such that 
$C \times_k k'$ is integral and its normalisation $C'$ of $C \times_k k'$ 
is $k'$-isomorphic to $\mathbb P^1_{k'}$: 
$$\varphi:C' \to C \times_k k' \to C.$$ 
In particular, $\varphi$ is a finite universal homeomorphism of degree $p$. 
By Proposition~\ref{p-cubic}(2)(3)(4), 
we can find a $k$-rational point $P$ around which $C$ is smooth over $k$. 
We set $D:=-P$ and let $\xi \in H^1(C, \MO_C(D))$ be a nonzero element 
whose existence is guaranteed by Serre duality. 
Since $C' \simeq \mathbb P^1_{k'}$ and $P':=\varphi^*P$ is a $k'$-rational point, 
we have that $H^1(C', \MO_{C'}(\varphi^*D))=0$ by Serre duality. 
Therefore, we can apply the construction as in Subsection~\ref{ss-construction} (cf. Notation \ref{n-general}). 
Then we obtain a cartesian diagram of regular projective $k$-varieties: 
$$\begin{CD}
S' @>\psi>> S\\
@VV\rho' V @VV\rho V\\
C' @>\varphi>> C.
\end{CD}$$
\end{nota}

\begin{lem}\label{l-C_2}
We use Notation~\ref{n-dP}. 
Then the following hold. 
\begin{enumerate}
\item If $p=2$, then $C_2$ is $k$-isomorphic to 
a conic curve in $\mathbb P^2_{k}$ or $\mathbb P^2_{k'}$. 
\item If $p=3$, then $C_2$ is $k'$-isomorphic to $\mathbb P^1_{k'}$. 
\end{enumerate}
\end{lem}

\begin{proof}
We have that 
$$(K_{S/C}+C_2) \cdot_k C_2=\deg_k(\omega_{C_2})-\deg_k(\rho^*\omega_{C}|_{C_2})
=\deg_k(\omega_{C_2}).$$
Moreover, it follows that 
$$\psi^*(K_{S/C}+C_2) \cdot_{k} C'_2=(K_{S'/C'}+pC'_2) \cdot_{k} C'_2$$
$$=(p-1)C'_2 \cdot_{k} C'_2=(p-1)\deg_{k}(\varphi^*D)=-p(p-1),$$
where the first equation holds by Lemma \ref{l-p-cover}, 
the third one by Lemma \ref{l-pm-sections} 
and the last one by $D=-P$. 
Since $\psi|_{C'_2}:C'_2 \to C_2$ is birational, it follows that 
$$\deg_k(\omega_{C_2})=
(K_{S/C}+C_2) \cdot_k C_2=\psi^*(K_{S/C}+C_2) \cdot_{k} C'_2=-p(p-1).$$
Since $C_2$ is a projective curve such that $\omega_{C_2}^{-1}$ is an ample invertible sheaf, 
it holds that $H^1(C_2, \MO_{C_2})=0$. 
Then, it follows from \cite[Lemma 10.6]{Kol13} that 
\begin{itemize}
\item $C_2$ is $K$-isomorphic to $\mathbb P^1_K$, or 
\item $C_2$ is $K$-isomorphic to a conic curve on $\mathbb P^2_K$, 
\end{itemize}
where $K=H^0(C_2, \MO_{C_2})$. 
In any case, it holds that $\deg_K(\omega_{C_2})=-2$. 
The natural morphisms $C' \simeq C'_2 \to C_2 \to C$ 
induce field extensions 
$$k=H^0(C, \MO_C) \subset H^0(C_2, \MO_{C_2}) \subset H^0(C'_2, \MO_{C'_2})=k',$$
which implies that $K=H^0(C_2, \MO_{C_2})$ is either $k$ or $k'$. 
Thus (1) holds.

We show (2). 
Since $p=3$, we get $\deg_k(\omega_{C_2})=-6$. 
Thus we have that 
$k'=H^0(C_2, \MO_{C_2})$ and $\deg_{k'}(\omega_{C_2})=-2$. 
Since $\psi|_{C'_2}:C'_2 \to C_2$ is birational and $\deg_{k'}(\omega_{C'_2})=\deg_{k'}(\omega_{C_2})$, 
$\psi|_{C'_2}$ is an isomorphism. 
Hence (2) holds. 
\end{proof}

\begin{thm}\label{t-bad-dP}
We use Notation~\ref{n-dP}. 
Then there exists an effective $\Q$-divisor $\Delta_S$ on $S$ 
such that $(S, \Delta_S)$ is klt and $-(K_S+\Delta_S)$ is ample. 
\end{thm}

\begin{proof}
By Proposition~\ref{p-criterion}, 
it suffices to show that $-D-K_C$ is ample and $(S, \frac{2}{p}C_2)$ is log canonical. 
The ampleness of $-D-K_C$ follows from $D=-P$ and $K_C \sim 0$. 
If $p=3$, then 
it follows from Lemma~\ref{l-C_2} that $(S, \frac{2}{3}C_2)$ is log canonical. 
If $p=2$, then 
$C_2$ is semi log canonical by Lemma~\ref{l-conic} and Lemma~\ref{l-C_2}. 
Therefore, $(S, C_2)$ is log canonical by inversion of adjunction 
(cf. \cite[Theorem 5.1]{Tanb}). 
\end{proof}

\begin{thm}\label{t-bad-dP2}
Let $k$ be an imperfect field whose characteristic $p$ is two or three. 
Then there exists a projective $\Q$-factorial klt surface $T$ over $k$ 
with $k=H^0(T, \MO_T)$ that satisfies the following properties. 
\begin{enumerate}
\item $-K_T$ is ample, 
\item $\rho(T)=1$, 
\item there is a Cartier divisor $M$ such that $M \not\sim 0$ and $pM \sim 0$, 
and 
\item there exists a finite universal homeomorphism 
$\mathbb P^2_{k'} \to T$, where $k \subset k'$ is 
a purely inseparable extension of degree $p$.  
\end{enumerate}
\end{thm}

\begin{proof}
We use Notation~\ref{n-dP}. 
There is a $\mathbb P^1$-bundle structure $\rho:S \to C$. 
Since $S'=\mathbb P_{\mathbb P^1_{k'}}(\MO_{\mathbb P^1} \oplus \MO_{\mathbb P^1}(1))$, 
we have the blow-down $f':S' \to \mathbb P^2_{k'}=:T'$ contracting $C'_2$. 
Thus, we get a commutative diagram 
$$\begin{CD}
S' @>\psi>> S\\
@VVf'V @VVfV\\
T' @>\psi_T>> T,
\end{CD}$$
where $\psi_T$ is a finite universal homeomorphism of degree $p$ and 
$f$ is the birational morphism to a projective normal surface $T$ 
satisfying $\Ex(f)=C_2$. 
Note that we can find such a surface $T$ as the Stein factorisation of 
$S \xrightarrow{\widetilde{\psi}} S' \xrightarrow{f'} T'$, where 
$\widetilde{\psi}$ is obtained by the factorisation $F^e_{S'}:S' \xrightarrow{\psi} S \xrightarrow{\widetilde{\psi}} S'$ of the $e$-th iterated absolute Frobenius morphism $F^e_{S'}$ for some $e$. 
Thus (4) holds. 
Since $f$ can be considered as the contraction of 
a $(K_S+\Delta_S)$-negative extremal ray, 
we have that $(T, f_*\Delta_S)$ is klt and $T$ is $\Q$-factorial 
(cf. \cite[Theorem 4.4]{Tanb}). 
In particular, $T$ is klt. 
By \cite[Theorem 4.4]{Tanb}, the assertions (1) and (2) hold. 
We get (3) by Proposition~\ref{p-cubic}(5)  and \cite[Theorem 4.4]{Tanb}. 
\end{proof}

\begin{rem}
The surface $T$ constructed in the proof of Theorem \ref{t-bad-dP2} 
has a unique singular point $t$. 
We proved that the singular point $t$ is klt in the proof of Theorem \ref{t-bad-dP2}. 
On the other hand, we see that $t$ is not a canonical singularity as follows. 
Let $a$ be the rational number defined by the following equation
\[
K_S=f^*K_T + a C_2. 
\]
It suffices to prove that $a<0$. 
The proof of Lemma \ref{l-C_2} implies that 
\[
(K_S+C_2) \cdot_k C_2= \deg_k(\omega_{C_2})=-p(p-1).
\]
On the other hand, it holds that 
\[
C_2 \cdot_k C_2=\deg_k(\MO_S(C_2)|_{C_2})=
\deg_k(\MO_{S'}(pC'_2)|_{C'_2})
\]
\[
=p \deg_k(\MO_{\mathbb P^1_{k'}}(-1))=-p^2,
\]
where the second equation holds by Lemma \ref{l-p-cover} 
and the third one follows from Lemma \ref{l-pm-sections} and Notation \ref{n-dP}. 
To summarise, we obtain $K_S \cdot_k C_2=p$ and $a=-1/p$, as desired. 
\end{rem}

\section{Pathological Mori fibre spaces}

\subsection{Mori fibre spaces with non-trivial torsion divisors}\label{ss-MFS1}

\begin{proof}(of Theorem~\ref{intro-patho2}) 
By Theorem~\ref{t-bad-dP2}, 
there exist a projective $\Q$-factorial klt surface $T$ over $k(t)$ 
with $H^0(T, \MO_T)=k(t)$ 
and a Cartier divisor $M$ on $T$ 
which satisfy the properties (1)--(4) in Theorem~\ref{t-bad-dP2}. 
We can find a projective morphism $f:V \to W$ and a Cartier divisor $D$ on $V$, 
where $W$ is a non-empty open subset $W$ of $\Spec\,k[t]$, 
$V \times_W \Spec\,k(t)=T$, and $D|_T=M$. 
In particular, the property (2) in the statement holds. 
After possibly shrinking $W$, 
we may assume that $V$ is normal, 
$f_*\MO_V=\MO_W$, 
$pD \sim 0$, and 
$-K_V$ is an $f$-ample $\Q$-Cartier $\Q$-divisor. 
Since $D|_T=M$, the property (4) in the statement holds. 
After possibly shrinking again, the property (1) (resp. (3)) in the statement holds 
by Lemma~\ref{l-rel-properties} (resp. Lemma~\ref{l-picard-one}). 
\end{proof}

\subsection{Mori fibre spaces with non-klt bases}\label{ss-MFS2}

The main purpose of this subsection is to show 
Theorem~\ref{t-imperfect-nonklt}, 
since it directly implies one of our main results: Theorem~\ref{intro-patho1}. 
In Part \ref{sss-setup}, we summarise notation. 
In Part \ref{sss-1st-step} and Part \ref{sss-2nd-step}, 
we run a suitable minimal model program that will be needed 
in the proof of Theorem~\ref{t-imperfect-nonklt}. 
In Part \ref{sss-final}, we prove Theorem~\ref{t-imperfect-nonklt} and 
Theorem~\ref{intro-patho1}.

\subsubsection{Setup}\label{sss-setup}

We use Notation~\ref{n-dP}. 
Assume that $[k:k^p]<\infty$. 
Let $M_C:=\MO_C(P)$, where $P$ is a $k$-rational point on $C$ 
around which $C$ is smooth over $k$. 
Let $M_S:=\rho^*M_C$, $M_{C'}:=\varphi^*M_C$, and 
$M_{S'}:=\psi^*\rho^*M_C$. 
We set  
$$X:=\mathbb P_S(\MO_S \oplus M_S), \quad R:=\mathbb P_C(\MO_C \oplus M_C)$$
$$X':=\mathbb P_{S'}(\MO_{S'} \oplus M_{S'}), \quad 
R':=\mathbb P_{C'}(\MO_{C'} \oplus M_{C'})$$
and obtain a cartesian diagram: 
$$\begin{CD}
X @>\rho_X >> R\\
@VV\pi V @VV\pi_R V\\
S @>\rho >> C,
\end{CD}$$
whose base change by $(-) \times_C C'$ can be written as  
$$\begin{CD}
X' @>\rho'_X >> R'\\
@VV\pi' V @VV\pi'_R V\\
S' @>\rho' >> C'. 
\end{CD}$$
Let $C^{\pm}$ be the sections of $\pi_R$ corresponding to 
the direct sum decomposition $\MO_C \oplus M_C$ 
such that $\MO_R(C^{\pm})|_{C^{\pm}}=\pm M_C$ if we identify $C$ with $C^{\pm}$. 
We set $S^{\pm}, C'^{\pm}$ and $S'^{\pm}$ to be the pull-backs of $C^{\pm}$ 
to $X, R'$ and $X'$, respectively.

Since $R' \simeq \mathbb P_{\mathbb P^1_{k'}}(\MO_{\mathbb P^1} \oplus \MO_{\mathbb P^1}(1))$, 
we have that $C'^-$ is a $(-1)$-curve on $R'$, 
i.e. $K_{R'} \cdot_{k'} C'^-=C'^- \cdot_{k'} C'^-=-1$. 
Let 
$$\theta':R' \to \mathbb P^2_{k'}=:Q'$$
be the blow-down contracting $C'^{-}$. 

Corresponding to $\theta'$, we can find a birational morphism 
$\theta:R \to Q$ to a projective surface $Q$ 
such that $\theta_*\MO_R=\MO_Q$ and $\Ex(\theta)=C^-$. 
Indeed, for a positive integer $e$ such that 
the $e$-th iterated absolute Frobenius morphism $F^e:R' \to R'=:R'^{[p^e]}$ 
factors through the induced morphism $R' \to R$, 
we define $Q$ as the normalisation of $Q'^{[p^e]}$ in $K(R)$, 
where $\theta'^{[p^e]}:R'^{[p^e]} \to Q'^{[p^e]}$ is defined as the same 
morphism as $\theta'$. 


Let $q:=\theta(\Ex(\theta))$ and $q':=\theta'(\Ex(\theta'))$.  
In the proof of Theorem~\ref{t-imperfect-nonklt}, we will run 
an $S^-$-MMP over $Q$ 
$$X=:X_0 \overset{f_0}\dashrightarrow X_1 \xrightarrow{f_1} X_2,$$
consisting of two steps: $f_0$ is a flip and $f_1$ is a divisorial contraction. 
To this end, we construct the corresponding $S'^-$-MMP 
in Part \ref{sss-1st-step} and Part \ref{sss-2nd-step}.

\subsubsection{The first step: flip}\label{sss-1st-step}

We use the same notation as in Part \ref{sss-setup}. 
Let $H_{X'}$ be an ample Cartier divisor on $X'$ and 
we define $\lambda'$ as 
$\lambda':=\sup\{\lambda \in \R_{\geq 0}\,|\, H_{X'}+\lambda S'^-\text{ is nef}\}$. 
Set 
$$L':=H_{X'}+\lambda' S'^-.$$
Since $S'^-$ is a smooth projective surface such that $-K_{S'^-}$ is ample, 
\begin{itemize}
\item any nef Cartier divisor $N$ on $S'^-$ is semi-ample 
(cf. \cite[Theorem 1.3]{Tan14}) and 
\item there are finitely many curves $\gamma_1, \cdots, \gamma_r$ on $S'$ 
such that ${\rm NE}(S')=\sum_{i=1}^r\R_{\geq 0}[\gamma_i]$ 
(cf. \cite[Theorem 2.14]{Tanb}). 
\end{itemize}
Hence, $\lambda'$ is a positive rational number and 
$L'|_{S'^-}$ is semi-ample. 
It follows from Keel's theorem (\cite[Proposition 1.6]{Kee99}) 
that $L'$ is semi-ample. 
Let 
$$g':X' \to Z'$$
be the birational contraction with $g'_*\MO_{X'}=\MO_{Z'}$ 
induced by $L'$. 
Since $L'|_{S'^-}$ is semi-ample and $(L'|_{S'^-}) \cdot (\pi'^*\rho'^*(c')|_{S'^-})>0$ 
for a closed point of $c' \in C'$, 
we have that $\Ex(g')$ is equal to the $(-1)$-curve $\Gamma'$ on $S'^-$. 
In particular, $g'$ is a small birational morphism.

We construct a flip of $g'$. 
Let 
$$h':Y' \to X'$$
be the blowup along $\Gamma'$. 
We have that $E':=\Ex(h')$ is isomorphic to 
$\mathbb P_{\Gamma'}(N_{\Gamma'/X'})$, 
where $N_{\Gamma'/X'}$ is the normal bundle, which is an extension of 
$N_{S'^-/X'}|_{\Gamma'}$ and $N_{\Gamma'/S'^-}$. 
Since 
$$S'^- \cdot_{k'} \Gamma'=-1, \quad (\Gamma' \text{ in } S'^-) \cdot_{k'} (\Gamma' \text{ in } S'^-)=-1,$$
the locally free sheaf $N_{\Gamma'/X'}$ is corresponding to an extension class 
$\alpha \in \text{Ext}^1_{\mathbb P^1}(\MO(-1), \MO(-1))=0$. 
Therefore, we get 
$N_{\Gamma'/X'} \simeq \MO_{\Gamma'}(-1) \oplus \MO_{\Gamma'}(-1)$. 
It follows that 
$E' \simeq \mathbb P_{k'}^1 \times_{k'} \mathbb P_{k'}^1$. 
Let $T'$ be the proper transform of $S'^-$ on $Y'$. 

\begin{lem}\label{l-422calculation}
The following hold. 
\begin{enumerate}
\item 
$K_{X'} \cdot \Gamma'=0$. 
\item 
$\MO_{Y'}(E')|_{E'} \simeq \MO_{\mathbb P^1 \times \mathbb P^1}(-1, -1)$ 
if we identify $E'$ with $\mathbb P_{k'}^1 \times_{k'} \mathbb P_{k'}^1$. 
\item 
$\MO_{Y'}(T')|_{E'} \simeq \MO_{\mathbb P^1 \times \mathbb P^1}(0, 1)$ 
if we identify $E'$ and $E' \to \Gamma'$ 
with $\mathbb P_{k'}^1 \times_{k'} \mathbb P_{k'}^1$ and 
its first projection, respectively. 
\end{enumerate}
\end{lem}

\begin{proof}
The assertion (1) follows from $S'^- \cdot_{k'} \Gamma'=-1$ and 
$$(K_{X'}+S'^-) \cdot_{k'} \Gamma'=K_{S'^-} \cdot_{k'} \Gamma'=-1,$$
where the latter equation follows from the fact that $\Gamma'$ is a $(-1)$-curve. 

We show (2). 
Since $K_{Y'}=h'^*K_{X'}+E'$, we have that 
{\small 
$$\MO_{\mathbb P^1 \times \mathbb P^1}(-2, -2) \simeq 
\MO_{Y'}(K_{Y'}+E')|_{E'}=\MO_{Y'}(h'^*K_{X'}+2E')|_{E'} \simeq \MO_{Y'}(2E')|_{E'},$$
}
where the last isomorphism holds because (1) implies 
$K_{X'}|_{\Gamma'} \sim 0$. 
Thus (2) holds.  

The assertion (3) follows from $h'^*S'^-=T'+E'$, 
(2) and 
$\MO_X(S'^-)|_{\Gamma'} \simeq  \MO_{\mathbb P^1}(-1)$. 
\end{proof}

\begin{lem}\label{l-rel-glob-nef}
Let $k_1$ be a field. 
Let $\zeta:Y_1 \to Z_1$ be a birational $k_1$-morphism 
of projective normal $k_1$-varieties. 
For any $\Q$-Cartier $\Q$-divisor $N$ on $Y_1$ 
and an ample $\Q$-Cartier $\Q$-divisor $H$ on $Z_1$, 
there exists a positive integer $m$ 
such that $(N+m\zeta^*H)|_{S}$ is big 
for any integral closed subscheme $S$ on $Y_1$ 
such that $S \not\subset \Ex(\zeta)$. 
In particular, if $N$ is $\zeta$-nef and $\zeta(\Ex(\zeta))$ is one point, then 
$N+m\zeta^*H$ is nef. 
\end{lem}

\begin{proof}
Let 
$$
I:=\{S\,|\,S\text{ is an integral closed subscheme of }Y_1\text{ such that  }
S\not\subset \Ex(\zeta)\}.
$$
If $S \in I$, the induced morphism $S \to \zeta(S)$ is birational. 
Therefore, for any $S\in I$, 
there is $n_S\in\mathbb Z_{>0}$ such that $(N+n_S\zeta^*H)|_{S}$ is big. 
Let $n_1:=n_{Y_1}$. 
By Kodaira's lemma, we may write 
$N+n_1\zeta^*H=A+D$ where $A$ is an ample $\Q$-Cartier $\Q$-divisor and $D$ is an effective $\Q$-divisor on $Y_1$. 
Let $D=\sum e_jD_j$ be the  decomposition into the irreducible components and we set $n_2:=\max_{D_j\in I}\{n_1, n_{D_j}\}$. 
For any $D_j\in I$, let $D_j^N$ be its normalisation and  we again apply Kodaira's lemma to $(N+n_2\zeta^*H)|_{D_j^N}$. 
Repeating the same procedure finitely many times, 
we can find $n\in\mathbb Z_{>0}$ such that 
$(N+n\zeta^*H)|_{S}$ is big for every $S\in I$. 
\end{proof}

For an ample Cartier divisor $H_{Z'}$ on $Z'$ and a sufficiently large integer $m$, 
we set 
$$M':=T'+mh'^*g'^*H_{Z'}.$$
We have that $M'$ is nef and big by Lemma \ref{l-422calculation}(3) and  Lemma~\ref{l-rel-glob-nef}. 
It follows from Lemma~\ref{l-rel-glob-nef} that $\mathbb E(M')=E'$, 
where we refer to \cite[Definition 0.1]{Kee99} 
for the definition of $\mathbb E(M')$. 
By \cite[Theorem 0.2]{Kee99}, we get the birational morphism
$$h'_1:Y' \to X'_1$$
with $(h'_1)_*\MO_{Y'}=\MO_{X'_1}$ induced by $M'$. 
By our construction, we have that $X'_1$ is $\Q$-factorial, 
$\rho(X'_1)=\rho(X'_0)=3$ (cf. \cite[Lemma 2.1]{CT}), 
$Y' \to Z'$ factors through $h'_1$, 
the fibre of the induced morphism $X'_1 \to Q'$ over $q'$ 
is set-theoretically equal to $S'^-_1 \cup \Gamma'_1$, 
and $\Gamma'_1 \not\subset S'^-_1$, where $\Gamma'_1:=h'_1(E')$ and 
$S'^-_1$ is the proper transform of $S'^-$. 
In particular, $S'^-_1$ is ample over $Z'$, hence $X'_1 \to Z'$ is 
a flip of $X' \to Z'$.

\subsubsection{The second step: divisorial contraction}\label{sss-2nd-step}

We use the same notation as in Part \ref{sss-setup} and Part \ref{sss-1st-step}.

\begin{lem}\label{l-normal}
The following hold. 
\begin{enumerate}
\item 
The normalisation of $S'^-_1$ is a universal homeomorphism. 
\item 
$-S'^-_1|_{S'^-_1}$ is ample. 
\end{enumerate}
\end{lem}

\begin{proof}
We show (1). 
Let $S'^-_{Y'}$ be the proper transform of $S'^-$ on $Y'$. 
Note that $\widetilde{h'}:S'^-_{Y'} \xrightarrow{\simeq} S'^-$ and 
the exceptional locus of $\widetilde{h'_1}:S'^-_{Y'} \to S'^-_1$ is equal to $\Gamma'_{Y'}$, 
where $\Gamma'_{Y'}:=(\widetilde{h'})^{-1}(\Gamma')$. 
Since $\Gamma'_{Y'} \simeq \Gamma' \simeq \mathbb P^1_{k'}$, 
we have that $\widetilde{h'_1}(\Gamma'_{Y'})$ is a $k'$-rational point and 
the induced morphism $\Gamma'_{Y'} \to \widetilde{h'_1}(\Gamma'_{Y'})$ 
is the same as the structure morphism $\mathbb P^1_{k'} \to \Spec\,k'$. 
In particular, any fibre of $\widetilde{h'_1}:S'^-_{Y'} \to S'^-_1$ is geometrically connected. 
Since $\widetilde{h'_1}$ factors through the normalisation 
$\nu_{S'^-_1}:(S'^-_1)^N \to S'^-_1$ of $S'^-_1$, 
any fibre of $\nu_{S'^-_1}$ is geometrically connected. 
Hence, (1) holds.

We show (2). 
Take a curve $B'$ on $Q'$ passing through $q'$ such that $|B'|$ is base point free. 
Then the inverse image $D'$ to $X'_1$ can be written as 
$$D'=aS'^-_1+F'$$
where $a>0$ and $F'$ is a nonzero effective $\Q$-divisor with $S'^-_1 \not\subset \Supp F'$. 
Take a general curve $G'$ on $S'^-_1$. 
Since $D' \cdot G'=0$ and $F' \cdot G'>0$, we have that $S'^-_1 \cdot G'<0$. 
Thus (2) holds by $\rho(S'^-_1)=1$.  
\end{proof}

Let $H_{X'_1}$ be an ample Cartier divisor on $X'_1$ and 
we define $\nu'$ by 
$\nu':=\sup\{\nu \in \R_{\geq 0}\,|\, H_{X'_1}+\nu S'_1\text{ is nef}\}.$ 
We have that $\nu'$ is a positive rational number and let 
$$N':=H_{X'_1}+\nu' S'^-_1.$$
Since we can find a positive integer $m$ 
such that $\MO_{X'_1}(mN')|_{S'^-_1} \simeq \MO_{S'^-_1}$ 
by Lemma~\ref{l-normal}, 
we have that $N'$ is semi-ample by Keel's theorem (\cite[Proposition 1.6]{Kee99}). 
Let 
$$f'_1:X'_1 \to X'_2$$
be the birational morphism induced by $N'$ with $(f'_1)_*\MO_{X'_1}=\MO_{X'_2}$. 
We also get a morphism $\alpha':X'_2 \to Q'$.

\subsubsection{Proof of Theorem~\ref{intro-patho1}}\label{sss-final}

\begin{thm}\label{t-imperfect-nonklt}
Let $k$ be an imperfect field  
whose characteristic $p$ is two or three. 
If $[k:k^p]<\infty$, then there exists a $k$-morphism 
$\alpha:X_2 \to Q$ 
of projective normal $k$-varieties, 
with $\alpha_*\MO_{X_2}=\MO_Q$ and $H^0(Q, \MO_Q)=k$, 
that satisfies the following properties: 
\begin{enumerate}
\item $X_2$ is a $\Q$-factorial threefold of Fano type, 
\item $Q$ is a projective $\Q$-factorial log canonical surface 
which is not klt, 
\item any fibre of $\alpha$ is geometrically irreducible of dimension one, a general fibre of $\alpha$ is $\mathbb P^1$, and  
\item $\rho(X_2/Q)=1$.  
\end{enumerate}
\end{thm}

\begin{proof}
We use the same notation as in Part \ref{sss-setup}, Part \ref{sss-1st-step} and 
Part \ref{sss-2nd-step}. 
We get the rational maps 
$$X=:X_0 \overset{f_0}\dashrightarrow X_1 \xrightarrow{f_1} X_2 \xrightarrow{\alpha} Q$$
corresponding to 
$$X'=:X'_0 \overset{f'_0}\dashrightarrow X'_1 \xrightarrow{f'_1}  X'_2 \xrightarrow{\alpha'}  Q',$$
where $X_0 \overset{f_0}\dashrightarrow X_1 \xrightarrow{f_1} X_2$ 
is an $S^-$-MMP over $Q$. 
Indeed, for a positive integer $e$ such that 
the $e$-th iterated absolute Frobenius morphism $F^e:X' \to X'=:X'^{[p^e]}$ 
factors through the induced morphism $X' \to X$, 
we define $X_i$ as the normalisation of $X'^{[p^e]}_i$ in $K(X)$, 
where $X'^{[p^e]} \dashrightarrow X'^{[p^e]}_i$ is the same birational map 
as $X' \dashrightarrow X'_i$. 

Since $X'_0$ is of Fano type, $f_0$ is small and $f_1$ is birational, 
we have that also $X_2$ is of Fano type by Lemma~\ref{l-fano-birat}. 
Thus (1) holds. 
Since $Q' \simeq \mathbb P^2_{k'}$ is $\Q$-factorial, so is $Q$ (cf. \cite[Lemma 2.5]{Tana}). 
We can write $K_R+bC^-=\theta^*K_Q$ for some $b \in \Q$. 
By $\theta^*K_Q \cdot_k C^-=0$, $(C^-)^2<0$ and 
\[
(K_R+C^-) \cdot_k C^- = \deg_k (\omega_{C^-})=0,
\]
we have that $b=1$. 
It follows from \cite[Corollary 2.13]{Kol13} that $Q$ is not klt but log canonical. 
Thus (2) holds. 
The assertion (3) follows from the construction, 
because the fibre of $X' \to Q'$ over $q'$ is an image of 
$\Ex(h') \simeq \mathbb P^1_{k'} \times_{k'} \mathbb P^1_{k'}$ and hence geometrically irreducible. 
Thanks to (3), we have that $\rho(X'_2/Q')=1$. 
The assertion (4) holds by $\rho(X'_2/Q')=1$ and Lemma~\ref{l-rel-inv}. 
\end{proof}

\begin{proof}(Proof of Theorem~\ref{intro-patho1})
We apply Theorem~\ref{t-imperfect-nonklt} for a field $k(t)$. 
Then there exists a $k(t)$-morphism 
$\alpha:X_2 \to Q$ 
of projective normal $k(t)$-varieties, 
with $\alpha_*\MO_{X_2}=\MO_Q$ and $H^0(Q, \MO_Q)=k(t)$,  satisfying the properties 
(1)--(4) in Theorem~\ref{t-imperfect-nonklt}. 
We can find projective $k$-morphisms 
$$V \xrightarrow{f} W \xrightarrow{g} T$$
of normal $k$-varieties such that 
$T$ is a non-empty open subset of $\Spec\,k[t]$ and 
$f \times_T \Spec\,k(t)=\alpha$. 
After possibly shrinking $T$, 
we may assume that 
\begin{itemize}
\item $V$ and $W$ are $\Q$-factorial by Lemma~\ref{l-rel-properties}, 
\item $V$ is klt by Lemma~\ref{l-rel-properties}, 
\item $W$ is not klt by Lemma~\ref{l-rel-properties}, and 
\item $f_*\MO_V=\MO_W$.  
\end{itemize} 
We set $W_1$ to be the subset of $W$ consisting of the points $w \in W$ 
such that $V_w$ is geometrically irreducible and of dimension one. 
By \cite[9.5.5 and 9.7.7]{EGAIV3}, $W_1$ is a constructible subset of $W$. 

\begin{claim}
There exists an open subset $W_2$ of $W$ such that 
$W_{\eta} \subset W_2 \subset W_1$, where $\eta$ is the generic point of $T$. 
\end{claim}

\begin{proof}(of Claim) 
For a subset $B$ of a set $A$, let $B^c:=A \setminus B$. 
By Theorem~\ref{t-imperfect-nonklt}(3), 
we have that $W_{\eta}  \subset W_1$. 
This inclusion implies that $\eta \not\in g(W_1^c)$, 
hence the constructible subset $g(W_1^c)$ of a curve $T$ 
is a proper closed subset of $T$. 
Thus the inclusions $W_{\eta} \subset W_2 \subset W_1$ hold 
for $W_2:=g^{-1}(g(W_1^c)^c)$. 
This completes the proof of Claim. 
\end{proof}

After replacing $W$ by $W_2$, 
we may assume that the fibre $V_w$ over any point $w \in W$ 
is geometrically irreducible and of dimension one. 
In particular, $\rho(V/W)=1$. 
This completes the proof of Theorem~\ref{intro-patho1}. 
\end{proof}


\end{document}